\documentclass[leqno,12pt]{article}
\usepackage{graphicx}
\usepackage{nips06,times}
\usepackage{amsmath}
\usepackage{amsthm}
\usepackage{amsfonts}
\usepackage{amssymb}
\usepackage{latexsym}
\usepackage{epsf}
\usepackage{amscd,amsxtra}
\usepackage{subfigure}
\usepackage[center]{titlesec}

\numberwithin{equation}{section} \numberwithin{figure}{section}
\def\<{\langle}
\def\>{\rangle}
\def\be{\begin{equation}}
\def\bes{\begin{equation*}}
\def\ee{\end{equation}}
\def\ees{\end{equation*}}
\def\bee{\begin{eqnarray}}
\def\bees{\begin{eqnarray*}}
\def\eee{\end{eqnarray}}
\def\eees{\end{eqnarray*}}
\def\m{\mathrm}
\def\diff{{\m d}}
\newtheorem{lemma}{\hspace{5mm}{\large L}{\small EMMA}}[section]%used for define "Lemma"
\newtheorem{thm}{\hspace{5mm}{\large T}{\small HEOREM}}[section]%used for define "Lemma"
\title{An Estimate of the Gap of the first two\\
Eigenvalues in the Schr\"odinger Operator}

\author{
Shing-Tung Yau\thanks{Research supported by the National Science Foundation and Chinese University of Hong Kong.} \medskip\\
Harvard University\\
Department of Mathematics\\
One Oxford Street\\
Cambridge, MA 02138, USA\\
\\
\\
\emph{Dedicated to Professor Louis Nirenberg on his 75th birthday}
}
\begin{document}

\maketitle

%\begin{abstract}
%Abstract here.
%\end{abstract}

\section*{Introduction}
\label{sec:Intro}

\hskip 3ex
In my previous paper \cite{2} with I. M. Singer, B. Wong and Stephen
Yau, I gave a lower estimate of the gap of the first 2 eigenvalues
of the Schr\"odinger operator in case the potential is convex. In
this note we note that the estimate can be improved if we assume the
potential is strongly convex. In particular if the Hessian of the
potential is bounded from below by a positive constant, the gap has
a lower bound independent of dimension. We also find gap when the
potential is not necessary convex.

\section{Convex potential}
\label{sec:1}

\hspace{5mm}Let $\lambda_1$ and $\lambda_2$ be the first and second eigenvalues
of the operator $\Delta - V$, and $u_1$ and $u_2$ be their
corresponding eigenfunctions:
\bee
\Delta u_1 - V u_1 = -\lambda_1 u_1, \\
\Delta u_2 - V u_2 = -\lambda_2 u_2. \nonumber \eee
It is well known that the first eigenfunction $u_1$ must be a
positive function (a theorem of Courant). On the other hand, the
second eigenfunction changes sign since $\int u_1u_2 = 0$. Therefore
$u_2$ changes sign.

One can estimate $\lambda_2-\lambda_1$ by the following formula:
\be%
\lambda_2-\lambda_1 = \inf_{\int fu_1=0} \frac{\int \left|\nabla f
\right|^2u_1^2}{\int f^2u_1^2},
\ee

Here, we take another approach to derive the estimate on
$\lambda_2-\lambda_1$.

Since $u_1>0$, $u=\frac{u_2}{u_1}$ is a well-defined smooth function
on $\Omega$. Using the Hopf lemma and the Malgrange preparation
theorem, one has the following

\

\begin{lemma}\label{lam:1}
$u=\frac{u_2}{u_1}$ is smooth up to the boundary. It satisfies the
Neumann condition on the boundary.
\end{lemma}

When (\ref{lam:1}) are Neumann problems, Lemma \ref{lam:1} is
trivial, when (\ref{lam:1}) are Dirichlet problems, we argue in the
following way.

Note $\frac{\partial}{\partial \nu}u_1|_{\partial\Omega}\neq 0$.
Therefore, by using the equation,
\bee%
\Delta u &=& \frac{\Delta u_2}{u_1} - \frac{u_2\Delta
u_1}{u_1^2}-2\nabla \ln u_1 \cdot
\nabla\Big(\frac{u_2}{u_1}\Big)\\
&=&\frac{u_1\Delta u_2 - u_2\Delta u_1}{u_1^2} - 2\nabla\ln u_1\cdot
\nabla\Big(\frac{u_2}{u_1}\Big)
\nonumber\\
&=&-(\lambda_2-\lambda_1)\frac{u_2}{u_1}-2\nabla\ln
u_1\cdot\nabla\Big(\frac{u_2}{u_1}\Big)
\nonumber\\
&=&-(\lambda_2-\lambda_1)u-2\nabla\ln u_1\cdot\nabla u.\nonumber
\eee

We have the Neumann boundary condition $\frac{\partial u}{\partial
\nu}|_{\partial \Omega}=0$. Let

\be \varphi_1=-\ln u_1 \ee %

so that

\be\label{equ:1.5} \Delta u = -(\lambda_2-\lambda_1)u +
2\nabla\varphi_1\cdot \nabla u.\ee

\begin{thm}
\label{thm:1}
  Suppose the Ricci curvature of $\Omega$ is nonnegative and $\partial
  \Omega$ is convex, and
  \be
  \left\{
  {\begin{array}{l}
    \Delta u = -(\lambda_2-\lambda_1)u + 2W\cdot\nabla u\\
    \frac{\partial}{\partial \nu}u|_{\partial \Omega}=0,
  \end{array}}
  \right.
  \ee
  where $W$ is a vector field such that
  \be\label{equ:1.7}
  W_{i,i}\geq\sqrt{\frac{c}{2}}>0
  \ee
  then
  \be\label{equ:1.8}
  \lambda_2-\lambda_1\geq\frac{\theta^2(\beta)}{\mathrm{diam}(\Omega)^2}+\beta\sqrt{c},
  \ee
  where
  $\theta(\beta)=\sin^{-1}\frac{1}{\sqrt{1+\frac{\beta}{\sqrt{2}-\beta}}}$
  and $0<\beta<\sqrt{2}$ arbitrary.
\end{thm}

  \textbf{{\large P}{\small ROOF}.} \hspace{1mm}Consider
  \be
  F=\left|\nabla u\right|^2+\alpha u^2 \textrm{ with } \alpha\geq 0.
  \ee
By computation, we have

  \be
  F_i=2u_ju_{ji}+2\alpha u u_i,
  \ee
  \bee
  \Delta F &=& F_{ii}=2u_{ji}u_{ji}+2u_ju_{jii}+2\alpha
  u_iu_i+2\alpha u u_{ii}\\
  &=& 2\left|\nabla\nabla u\right|^2 + 2\nabla u \cdot \nabla\Delta u
  +\sum_{i,j}R_{ij}u_iu_j+2\alpha\left|\nabla u\right|^2+2\alpha
  u\Delta u
  \nonumber\\
  &=& 2\left|\nabla\nabla u\right|^2+2\nabla u\cdot\nabla(-(\lambda_2-
  \lambda_1)u-2W\cdot\nabla u)+\sum_{i,j}R_{ij}u_iu_j
  \nonumber\\
  &&+2\alpha\left|\nabla u\right|^2+ 2\alpha u(-(\lambda_2-\lambda_1)u-2W\cdot\nabla u)
  \nonumber\\
  &=&2\left|\nabla\nabla
  u\right|^2+\sum_{i,j}R_{ij}u_iu_j-2\Big((\lambda_2-\lambda_1)\left|\nabla
  u\right|^2
  \nonumber\\
  &&+\sum_{i,j}(W_{i,j}+W_{j,i})u_iu_j+2\sum_{i,j}W_iu_{ij}u_j\Big)
  \nonumber\\
  &&+2\alpha \left|\nabla
  u\right|^2-2\alpha((\lambda_2-\lambda_1)u^2+2u\nabla W\cdot\nabla
  u)
  \nonumber\\
  &=&2\left|\nabla\nabla
  u\right|^2+\sum_{ij}R_{ij}u_iu_j-2(\lambda_2-\lambda_1-\alpha)\left|\nabla
  u\right|^2
  \nonumber\\
  &&+2\sum_{i,j}(W_{i,j}+W_{j,i})u_iu_j-2\alpha(\lambda_2-\lambda_1)u^2+2W\cdot\nabla
  F.
  \nonumber
  \eee
  If $R_{ij}\geq 0$ and

  \be
  W_{i,i}\geq \sqrt{\frac{c}{2}},
  \ee
  then
  \bee
  \Delta F-2W\cdot\nabla F&\geq& 2\left|\nabla\nabla
  u\right|^2-2\Big(\lambda_2-\lambda_1-\alpha-2 \sqrt{\frac{c}{2}}
  \Big)\left|\nabla u\right|^2\\
  & & - 2 \alpha (\lambda_2 - \lambda_1)u^2. 
  \nonumber
  \eee

  First, we need to derive a universal lower bound for
  $\lambda_2-\lambda_1$.

\noindent(1) Let $\alpha=0$.

    If $F$ attains the maximum at the boundary point, say $x_0$,
    then $\frac{\partial}{\partial \nu}F(x_0)\geq 0$.

    Take a local orthonormal frame $(e_1,\ldots,e_n)$ near $x_0$
    such that $\nu=e_n$. From the definition of Hessian and second
    fundamental form, we have

    \bee
      u_{in} &=& e_ie_nu-(\nabla_{e_i}e_n)u\\
      &=& -(\nabla_{e_i}e_n)u  \hspace{1mm}\textrm{ since } \hspace{1mm}u_{\nu}=0 \nonumber\\
      &=& -\sum^{n-1}_{j=1}h_{ij}u_i.\nonumber
    \eee
    \bee
      F_{\nu} &=& 2\sum_j u_ju_{j\nu}\\
      &=& -2\sum h_{ij}u_iu_j \nonumber\\
      &\leq& 0 \hspace{1mm} \textrm{ \ by the convexity of }\hspace{1mm} \partial \Omega.
      \nonumber
    \eee

    This implies that $u_1=\ldots=u_{n-1}=0$, hence $\nabla u=0$ at
    $x_0$. Therefore, we have

    \be
    F\equiv 0.\nonumber
    \ee

    Thus $u$ is a constant which is impossible.

    If $F$ attains the maximum at the interior point, say $x_0$,
    then $\nabla u(x_0)\neq 0$. Otherwise, we have the same
    conclusion as above.

    At $x_0$,
    \bee
    0 &\geq& \Delta F(x_0)\\
    &\geq& 2\left|\nabla\nabla
    u\right|^2-2(\lambda_2-\lambda_1)\left|\nabla u\right|^2
    \nonumber\\
    &&+4\sqrt{\frac{c}{2}}\left|\nabla u\right|^2 \hspace{1mm}\textrm{ since }\hspace{1mm}
    \nabla F(x_0) = 0.\nonumber
    \eee
    The last inequality is equivalent to the following:
    \be
    \Big((\lambda_2-\lambda_1)-2\sqrt{\frac{c}{2}}\Big)\left|\nabla
    u\right|^2 \geq |\nabla\nabla u|^2 \geq 0,
    \ee
    which says that
    \be
    (\lambda_2-\lambda_1)\geq \sqrt{2c} \hspace{1mm}\textrm{ since } \hspace{1mm}\nabla
    u(x_0) \neq 0.
    \ee

    \noindent (2) Now, take $\alpha=\lambda_2-\lambda_1-\beta\sqrt{c}>0$

    From the universal lower bound, we can take
    $\beta=\sqrt{2}-\delta$ for any small $\delta>0$ in the
    following argument.

    \noindent\textbf{Case 1.}\hspace{2mm} If $x_0\in \partial \Omega$, then $\frac{\partial}{\partial \nu}
    F(x_0)\geq 0$.

    \bee
    F_v&=&2\sum_j u_j u_{j\nu} + 2\alpha u u_{\nu}\\
    &=&-2\sum h_{ij}u_iu_i \nonumber\\
    &\leq&0 \hspace{1mm}\textrm{ by the convexity of } \partial \Omega.\nonumber
    \eee
    This implies that $u_1=\ldots=u_{n-1}=0$, hence $\nabla u=0$ at
    $x_0$. Therefore, we have
    
    \be
    F\leq \sup \alpha u^2.
    \ee

\noindent    \textbf{Case 2.} \hspace{2mm}$x_0\in \mathop{\Omega}\limits^\circ$ and

    \be
    \textrm{(a)}\quad \nabla u(x_0)=0.
    \ee
    Then by the definition
    \be
    F(x_0)=\left|\nabla u\right|^2(x_0)+\alpha u^2(x_0)=\alpha
    u^2(x_0)\leq \alpha \sup u^2.
    \ee
    Hence
    \be
    \left|\nabla u\right|^2+\alpha u^2=F\leq \alpha \sup u^2.
    \ee

   \noindent \textbf{Case 3.} \hspace{2mm}$x_0 \in \mathop{\Omega}\limits^\circ$ and
    \be
    \textrm{(b)}\quad \nabla u(x_0)\neq 0.
    \ee
    Using
    \be
    0=F_i(x_0)=2u_ju_{ji}+\alpha u u_i = 2u_j(u_{ij}+\alpha u
    g_{ij})
    \ee
    and rotating normal coordinates centered at $x_0$, we may assume
    \bee
      &&u_1(x_0)\neq 0,\\
      &&u_i(x_0)=0,\hspace{1mm} i=2,\ldots,n.\nonumber
    \eee
    Then
    \be
    u_{11}+\alpha u=0
    \ee
    which implies
    \be
    u_{11}^2=\alpha^2u^2
    \ee
    so that
    \be
    \sum u_{ij}^2 \geq \alpha^2 u^2.
    \ee
    Hence
    \bee
    0&\geq&2\left|\nabla\nabla
    u\right|^2-2(\lambda_2-\lambda_1-\alpha) \left|\nabla u\right|^2
    - 2\alpha(\lambda_2-\lambda_1)u^2+4\sqrt{\frac{c}{2}} \left|\nabla
    u\right|^2\\
    &\geq&-2(\lambda_2-\lambda_1-\alpha-\sqrt{2c}) \left|\nabla
    u\right|^2 - 2\alpha(\lambda_2-\lambda_1-\alpha)u^2 \nonumber
    \eee
    Then
    \bee
    0&\geq&-2(\lambda_2-\lambda_1-\alpha-\sqrt{2c}) \left|\nabla
    u\right|^2 - 2\alpha(\lambda_2-\lambda_1-\alpha)u^2 \\
    &=&2(-\beta\sqrt{c}+\sqrt{2c})\left|\nabla
    u\right|^2 - 2\alpha\beta\sqrt{c}u^2,\nonumber
    \eee
    which implies
    \be
    (-\beta+\sqrt{2})\left|\nabla
    u\right|^2-\alpha\beta u^2 \leq 0
    \ee
    and if $\beta<\sqrt{2}$, at $x_0$
    \be
    \left|\nabla
    u\right|^2 \leq \frac{\alpha\beta}{-\beta +\sqrt{2}}u^2.
    \ee
    Hence, if $\beta<\sqrt{2}$, then at $x_0$
    \be
    F=\left|\nabla
    u\right|^2+\alpha
    u^2\leq\alpha\Big(1+\frac{\beta}{\sqrt{2}-\beta}\Big)u^2
    \ee
    so that
    \be
    F=\left|\nabla
    u\right|^2+\alpha
    u^2\leq\alpha\Big(1+\frac{\beta}{\sqrt{2}-\beta}\Big)\sup u^2,
    \ee
    which covers all the cases.

  Hence
  \be
  \frac{\left|\nabla u\right|}
  {\sqrt{\alpha(1+\frac{\beta}{\sqrt{2}-\beta})\sup u^2-\alpha
  u^2}}\leq 1.
  \ee
  Normalizing so that $\sup u^2=1$ and integrating along a shortest
  straight line $\gamma$ from $x_1$ where
  $\left|u(x_1)\right|=\sup\left|u\right|$ to the nodal set
  $\{u=0\}$, we obtain
  \bee
  \textrm{diam}(M)&\geq&\int_{\gamma}
  \frac{\left|\nabla u\right|}
  {\sqrt{\alpha(1+\frac{\beta}{\sqrt{2}-\beta})-\alpha
  u^2}}\\
  &\geq&\frac{1}{\sqrt{\alpha}}
  \int_0^1\frac{du}{\sqrt{1+\frac{\beta}{\sqrt{2}-\beta}-u^2}}
  \nonumber\\
  &=&\frac{1}{\sqrt{\alpha}}\sin^{-1}
  \frac{1}{\sqrt{1+\frac{\beta}{\sqrt{2}-\beta}}}
  \nonumber
  \eee
  so that
  \be
  \lambda_2-\lambda_1-\beta\sqrt{c}=\alpha\geq\Bigg(
  \sin^{-1}\frac{1}{\sqrt{1+\frac{\beta}{\sqrt{2}-\beta}}}\Bigg)^2
  \frac{1}{\textrm{diam}(M)^2},
  \ee
  \be
  \lambda_2-\lambda_1\geq\frac{\theta^2(\beta)}{\textrm{diam}(M)^2}+\beta
  \sqrt{c},
  \ee
  where
  \be
  \theta(\beta)=\sin^{-1}\frac{1}{{\sqrt{1+\frac{\beta}{\sqrt{2}-\beta}}}}
  \ee

  This finishes the proof of Theorem \ref{thm:1}. 
  
   Formula (\ref{equ:1.5}) will satisfy the hypothesis of Theorem
  \ref{thm:1} if the Hessian of $\varphi_1$ has a lower bound
  (\ref{equ:1.7}). This will be proved in section two for convex
  domain.

\

  \begin{thm}
    \label{thm:2}
    For a convex domain $\Omega$ with a potential $V$ whose Hessian
    has a lower bound $c>0$. Then (\ref{equ:1.8}) holds.
  \end{thm}
  
  \

  \section{Nonconvex Potential}
  
\hspace{5mm}  For the first eigenfunction $u_1$ defined on the domain $\Omega$
  in Euclidean space, we know that the Hessian of $\varphi=-\log
  u_1$ tends to infinity if $\partial \Omega$ is strictly convex and
  $u_1=0$ on $\partial\Omega$. Since
  
  \be\label{equ:2.1}
  \Delta\varphi=\left|\nabla\varphi\right|^2-V+\lambda_1,
  \ee
  we deduce
  \be
  \Delta\frac{\partial^2\varphi}{\partial x_i^2} =
  2\sum\Big(\frac{\partial^2\varphi}{\partial x_i\partial x_j}
  \Big)^2+2\nabla\varphi\cdot\nabla\Big(\frac{\partial^2\varphi}
  {\partial x_i^2}\Big)-\frac{\partial^2 V}{\partial x_i^2}.
  \ee

  If $\frac{\partial^2V}{\partial x_i^2}\geq c>0$ in $\Omega$, then
  we can argue from (\ref{equ:2.1}) that at point $x\in\Omega$ where
  $\frac{\partial^2\varphi}{\partial x_i^2}$ is minimum,
  $\frac{\partial^2\varphi}{\partial x_i\partial x_j}=0$ for $j\neq
  i$ and
  \be
  \label{equ:2.3}
  2\Big(\min_i\frac{\partial^2\varphi}{\partial x_i^2}\Big)^2 \geq
  \frac{\partial^2V}{\partial x_i^2} \geq c>0.
  \ee

  The continuity argument here was used by me in 1980 to handle the
  log concavity of $u_1$. By looking at $tV+\frac{(1-t)c\sum
  x_i^2}{2n}$, we know that when $t=0$, $\min_i
  \frac{\partial^2\varphi}{\partial x_c^2}\geq\sqrt{\frac{c}{2}}>0$.
  It follows from (\ref{equ:2.3}) that this must be valid when $t=1$
  also.

\

  \begin{thm}
    \label{thm:2.1}
    For a Dirichlet problem with $\frac{\partial^2 V}{\partial
    x_i^2} \geq c > 0$, the first eigenfunction $u_1$ satisfies the
    inequality $-\frac{\partial^2\log u_1}{\partial x_i^2}\geq
    \sqrt{\frac{c}{2}}>0$.
  \end{thm}

  We shall now treat the case when $V$ is not necessary convex. We
  shall assume Neumann boundary condition.

  First of all, we give an upper bound for for $\Delta\varphi$. From
  (\ref{equ:2.1}), it is trivial to verify that
  
  \be
  \Delta(\Delta\varphi)=2\nabla\varphi\cdot\nabla(\Delta \varphi)+2
  \left|\nabla\nabla\varphi\right|^2-\Delta V.
  \ee

  Since $\left|\nabla\nabla\varphi\right|^2\geq \frac{1}{n}
  (\Delta \varphi)^2$, we conclude that if $\Delta \varphi$ achieves
  its maximum in the interior of $\Omega$,
  
  \be
  (\Delta\varphi)^2\leq\frac{n\sup\Delta V}{2}.
  \ee

  On the other hand, if $\Delta\varphi$ achieves its maximum on the
  boundary $\partial\Omega$,
  \be
  \frac{\partial(\Delta\varphi)}{\partial \nu}\leq 0.
  \ee

  From (\ref{equ:2.1}) and that $\frac{\partial \varphi}{\partial
  \nu}=0$, we conclude that
  \be\label{equ:2.7}
  \sum_{i\neq \nu}\varphi_i\varphi_{i\nu}\leq\frac{\partial V}{\partial
  \nu}.
  \ee

  If the second fundamental of $\partial \Omega$ has eigenvalue
  greater than $\lambda>0$, we conclude from (\ref{equ:2.7}) that
  \be
  \left|\nabla\varphi\right|^2\leq\frac{1}{\lambda}\frac{\partial
  V}{\partial \nu}.
  \ee

  Therefore
  \bee
  \nabla\varphi&=&\left|\nabla \varphi\right|^2-V+\lambda_1\\
  &\leq&\frac{1}{\lambda}\frac{\partial V}{\partial \nu} - V + \lambda_1.
  \nonumber
  \eee

  \begin{thm}
    \label{thm:2.2}
    For the Neumann problem on a convex domain $\Omega$ whose
    boundary have principle curvature greater than $\lambda>0$. Then
    either
    \bee
    &&\nabla \varphi \leq\frac{n}{2}\sqrt{\sup_{\Omega}\Delta V}
    \textrm{\quad or}
    \nonumber\\
    &&\Delta\varphi \leq \sup_{\partial \Omega}\Big(\frac{1}{\lambda}\frac{\partial V}{\partial
    v}-V\Big)+\lambda_1.
    \nonumber
    \eee
    In particular for $\varphi=-\log u_1$, $\left|\nabla
    \varphi\right|^2-V +\lambda_1\leq \frac{n}{2}\sqrt{\sup_{\Omega}\nabla
    V}$ or $\sup\big(\frac{1}{\lambda}\frac{\partial V}{\partial
    v}\big)+\lambda_1$.
  \end{thm}

  In order to obtain lower estimate of the Hessian of $\varphi$, we
  argue as follows.

  For simplicity we shall assume that our domain is the ball in
  $R^n$. We shall use polar coordinate so that
  \be
  \Delta = \frac{\partial^2}{\partial
  r^2}+\frac{n-1}{r}\frac{\partial}{\partial
  r}+\frac{1}{r^2}\Delta_{\theta}.
  \ee

  Therefore the operator $\Delta_{\theta}$ commutes with $\Delta$ and
  we obtain
  \bee
  \label{equ:2.11}
  \Delta(\Delta_{\theta}\varphi)&=&2\varphi_r(\Delta_{\theta}\varphi)_r+
  2r^{-2}\sum_i
  \varphi_{\theta_i}(\Delta_{\theta}\varphi)_{\theta_i}\\
  &&+2(n-2)r^{-2}\sum\varphi_{\theta_i}^2+2\sum\varphi_{r\theta_i}^2
  \nonumber\\
  &&+2r^{-2}\sum\varphi_{\theta_i\theta_j}^2-\Delta_{\theta}V.\nonumber
  \eee

  Since we assume the Neumann boundary condition, $\varphi_r=0$
  along the boundary and so $(\Delta_{\theta}\varphi)_r=0$ along the boundary.
  By the sharp maximum principle, we can assume that $\Delta_{\theta}\varphi$
  achieves its
  maximum in the interior of $\Omega$ which implies by
  (\ref{equ:2.11}) that
  \be
  \sup\Delta_{\theta}\varphi\leq\frac{(n-1)^{1/2}}{\sqrt{2}}r\sup_{\Omega}
  (\Delta_{\theta}V)^{1/2}_+.
  \ee

  If we compute the upper bound of the spherical Hessian of
  $\varphi$, we can apply the same argument to find
  \be
  \label{equ:2.13}
  \sup_{\Omega}\frac{\partial^2\varphi}{\partial \theta_i^2} \leq
  \frac{1}{8}+r\sup_{\Omega}\Big(\frac{r\partial^2 V}{\partial \theta_i^2}\Big)^{1/2}_+.
  \ee

  In order to obtain estimate of the full Hessian of $\varphi$, we
  use the equation
  \bee
  \label{equ:2.14}
  \Delta\Big(\frac{r\partial \varphi}{\partial
  r}\Big)&=&2\Delta\varphi+\frac{r\partial(\Delta \varphi)}{\partial r}
  \\
  &=&2\Delta\varphi+2\nabla\varphi\cdot\nabla\Big(\frac{r
  \partial\varphi}{\partial r}\Big) - 2\left|\nabla\varphi\right|^2
  -\frac{r\partial V}{\partial r}\nonumber\\
  &=&-2V+2\lambda_1-\frac{r\partial V}{\partial
  r}+2\nabla\varphi\cdot \nabla\Big(\frac{r\partial\varphi}{\partial
  r}\Big).\nonumber
  \eee

  Similarly
  
  \bee
  \Delta\Big[r\frac{\partial}{\partial r}
  \Big(r\frac{\partial\varphi}{\partial r}\Big) \Big]&=& 2\Delta\Big(
  \frac{r\partial\varphi}{\partial r}\Big)+r\frac{\partial}{\partial
  r}\Delta\Big(r\frac{\partial \varphi}{\partial r}\Big)\\
  &=&-4V+4\lambda_1-2r\frac{\partial V}{\partial r}+
  4\nabla\varphi\cdot \nabla\Big(\frac{r\partial\varphi}{\partial r}\Big)
  \nonumber\\
  &&-2r\frac{\partial}{\partial r}(V-\lambda_1)-r\frac{\partial}{\partial
  r} \Big(r\frac{\partial V}{\partial r}\Big)\nonumber\\
  &&+2\Big|\nabla\Big(r\frac{\partial u}{\partial r}\Big)\Big|^2+
  2\nabla\varphi\cdot\nabla\Big(r\frac{\partial}{\partial
  r}\Big(r\frac{\partial u}{\partial r}\Big)\Big)\nonumber\\
  &&-4\nabla\varphi\cdot\nabla\Big(r\frac{\partial\varphi}{\partial
  r}\Big).\nonumber
  \eee

  Hence of $r\frac{\partial}{\partial r}\Big(r\frac{\partial\varphi}{\partial
  r}\Big)$ achieves its maximum in the interior of $\Omega$,
  
  \be
  2\Big|\nabla\Big(r\frac{\partial\varphi}{\partial r}\Big)\Big|^2
  \leq r\frac{\partial}{\partial r}\Big(r\frac{\partial V}{\partial
  r}\Big)+4r\frac{\partial V}{\partial r}+4V-4\lambda_1.
  \ee

  Hence in this case,
  
  \be
  \label{equ:2.17}
  \sup r\frac{\partial}{\partial r}\Big(r\frac{\partial\varphi}{\partial
  r}\Big) \leq \sup_{\Omega}\sqrt{\Big(\frac{1}{2}r\frac{\partial}{\partial r}
  \Big(\frac{\partial V}{\partial r}\Big)+2r\frac{\partial V}{\partial r} +
  2V-\lambda_1\Big)}.
  \ee

  If $r\frac{\partial}{\partial r}\Big(r\frac{\partial\varphi}{\partial
  r}\Big)$ achieves its maximum at the boundary of $\Omega$, we note
  that
  
  \bee
  r\frac{\partial}{\partial r}\Big(r\frac{\partial\varphi}{\partial r}
  \Big)&=&\frac{\diff^2\varphi}{\diff r^2}+\frac{n-1}{r}\frac{\partial\varphi}
  {\partial r}-\frac{n-2}{r}\frac{\partial\varphi}{\partial r}\\
  &=&\Delta\varphi-\frac{1}{r^2}\Delta_{\theta}\varphi-
  \frac{n-2}{r}\frac{\partial\varphi}{\partial r}\nonumber\\
  &=&\left|\nabla\varphi\right|^2-V+\lambda_1- \frac{1}{r^2}
  \Delta_{\theta}\varphi-\frac{n-2}{r}\frac{\partial\varphi}{\partial
  r}.\nonumber
  \eee

  Since $\frac{\partial\varphi}{\partial r}=0$ along the boundary
  and $\frac{\partial}{\partial r}\Big(r\frac{\partial}{\partial r}\Big(r
  \frac{\partial\varphi}{\partial r}\Big)\Big)\geq 0$ at the maximum
  point,
  
  \bee
  0&\leq&-\frac{2}{r^3}\left|\nabla_{\theta}\varphi\right|^2
  -\frac{\partial V}{\partial r}+
  \frac{2}{r^3}\Delta_{\theta}\varphi -\frac{n-2}{r}\frac{\partial^2\varphi}
  {\partial r^2}\\
  &=&-\frac{2}{r^3}\left|\nabla_{\theta}\varphi\right|^2
  -\frac{\partial V}{\partial r}+\frac{2}{r}\Big(\Delta\varphi
  -\frac{\partial^2\varphi}{\partial r^2}\Big)-
  \frac{n-2}{r}\frac{\partial^2\varphi}
  {\partial r^2}\nonumber\\
  &=&-\frac{2}{r^3}\left|\nabla_{\theta}\varphi\right|^2
  -\frac{\partial V}{\partial r}+\frac{2}{r}\Big(
  \frac{1}{r^2}\left|\nabla_{\theta}\varphi\right|^2-V+\lambda_1\Big)-
  \frac{n}{r}\frac{\partial^2\varphi}{\partial r^2}.\nonumber
  \eee

  Hence in this case
    \be
  \label{equ:2.20}
  \sup r\frac{\partial}{\partial r}\Big(r\frac{\partial\varphi}{\partial
  r}\Big) \leq \frac{1}{n}\sup_{\partial \Omega}\Big[-r^3\frac{\partial V}{\partial
  r}-2r^2(V-\lambda_1)\Big].
  \ee

  Hence either (\ref{equ:2.17}) or (\ref{equ:2.20}) hold.

  Note that since $\Delta\varphi$ is the sum of the Hessian of
  $\varphi$ in radial and spherical directions and sum we have upper
  estimate of Hessian in these directions, we have also lower
  estimate of them in terms of $\Delta\varphi$.

\

  \begin{thm}
    For the Neumann problem when $\Omega$ is a ball, and $\varphi=-\log
    u_1$, (\ref{equ:2.13}) holds for spherical Hessian and either
    (\ref{equ:2.17}) or (\ref{equ:2.20}) hold for radial Hessian.
  \end{thm}

  To obtain the full Hessian estimate of $\varphi$, we need to
  control $\varphi_{r\theta}$ and then can be accomplished as
  follows:

  Call $\psi=r\frac{\partial\varphi}{\partial r}$. Then according to
  equation (\ref{equ:2.14}), we compute
  
  \bee
  \label{equ:2.21}
  \Delta(\left|\nabla\psi\right|^2+c\psi^2)
  &=&
  2\sum\psi^2_{ij}+2\nabla\psi\nabla(\Delta\psi)+2c\left|\nabla\psi\right|^2
  +2c\psi\Delta\psi\\
  &=&2\sum\psi^2_{ij}-4\nabla\psi\cdot\nabla
  V-2\nabla\psi\cdot\nabla\Big(r\frac{\partial V}{\partial
  r}\Big)\nonumber\\
  &&+4\sum\varphi_i\psi_{ij}\psi_j+4\sum\psi_i\varphi_{ij}\psi_j
  \nonumber\\
  &&+2c\left|\nabla\psi\right|^2+2c\Big(-2V+2\lambda_1-r\frac{\partial
  V}{\partial r}\Big)\psi \nonumber\\
  &&+4c\psi\nabla\varphi\nabla\psi.\nonumber
  \eee

  If $\sup(\left|\nabla\psi\right|^2+c\psi^2)$ occurs in the
  interior, we obtain from (\ref{equ:2.21})
  
  \bee
  \label{equ:2.22}
  0&\geq&2\sum\psi^2_{ij}-4\nabla\psi\cdot\nabla
  V-2\nabla\psi\cdot\nabla\Big(r\frac{\partial V}{\partial r}\Big)\\
  &&+4\sum\psi_i\varphi_{ij}\psi_j+2c\left|\psi\right|^2\nonumber\\
  &&-4cV\psi+4c\lambda_1\psi-2cr\frac{\partial V}{\partial
  r}\psi.\nonumber
  \eee

  Note that
  \be
  \label{equ:2.23}
  \sum\psi_i\varphi_{ij}\psi_j=\psi^2\varphi_{rr}+2\psi_r\sum\varphi_{r\theta_j}\psi_{\theta_j}
  + 2\sum\psi_{\theta_i}\varphi_{\theta_i\theta_j}\psi_{\theta_j}.
  \ee

  Since we have already estimate $\varphi_{rr}$, $\psi_r$ and
  $\varphi_{\theta_i\theta_j}$, we conclude that
  $\sum\varphi_i\varphi_{ij}\varphi_j$ can be estimated by
  $\left|\nabla\psi\right|^2$. By choosing $C$ large enough, we
  conclude from (\ref{equ:2.23}) $\left|\nabla\psi\right|^2+c\psi^2$
  can be estimated from the information of $V$, $\nabla V$ and
  $\nabla\nabla V$.

  If $\left|\nabla\psi\right|^2+c\psi^2$ achieves its maximum on the
  boundary of $\Omega$,
  \be
  0\leq 2\sum\psi_j\psi_{j\nu}+2\psi\psi_{\nu}.
  \ee

  Note $\psi=0$ on $\partial\Omega$, and hence
  \bee
  \label{equ:2.25}
  0&\leq&\sum_j\psi_j\psi_{j\nu}\\
  &=&\psi_{\nu}\psi_{\nu\nu}\nonumber\\
  &=&\psi_{\nu}(\Delta\psi)-H\psi^2_{\nu}\nonumber\\
  &=&\psi_{nu}\Big(-2V+2\lambda_1-r\frac{\partial V}{\partial r}\Big) +
  2\varphi_v\psi_{\nu}^2-H\psi_{\nu}^2\nonumber
  \eee
  where $H$ is the mean curvature of $\partial\Omega$.

  As $\varphi_{\nu}=0$ on $\partial\Omega$, we conclude that if $\left|\nabla
  \psi\right|^2+c\psi^2$ achieves its maximum on $\partial\Omega$,
  
  \be
  \psi_{\nu}^2+c\psi^2\leq\sup_{\partial\Omega}\frac{1}{H^2}\Big(-2V+2\lambda_1-r\frac{\partial
  V}{\partial r}\Big)^2
  \ee
  
  \

  \begin{thm}
    If $\psi=r\frac{\partial\varphi}{\partial r}$,
    $\left|\nabla\varphi\right|$ can be estimated by $V$, $\nabla V$,
    $\nabla\nabla V$ using (\ref{equ:2.22}), (\ref{equ:2.23}) and
    (\ref{equ:2.25}).
  \end{thm}

  This completes estimates for the full Hessian of $\varphi$.

  Incidently (\ref{equ:2.14}) shows that
  
  \be
  \Delta\Big(r\frac{\partial \varphi}{\partial r}-2\varphi\Big)
  =2\nabla\varphi\cdot\nabla\Big(r\frac{\partial \varphi}{\partial
  r}-2\varphi\Big)+2\left|\nabla\varphi\right|^2-r\frac{\partial
  V}{\partial r}.
  \ee

  Suppose we want to find an upper estimate of
  $r\frac{\partial\varphi}{\partial r}-2\varphi$, we can proceed as
  follows. For any function $f$ such that
  
  \be
  \Delta f-\frac{1}{2}\left|\nabla f\right|^2-r\frac{\partial V}{\partial
  r}\geq 0
  \ee
  we find that at an interior maximum point of $r\frac{\partial\varphi}{\partial r}-2\varphi
  +f$, we have
  
  \bee
  \label{equ:2.29}
  0&\geq& 2\left|\nabla\varphi\right|^2-2\nabla\varphi\cdot\nabla
  f-r\frac{\partial V}{\partial r} + \Delta f\\
  &=&2\Big|\nabla\varphi-\frac{1}{2}\nabla
  f\Big|^2-\frac{1}{2}\left|\nabla f\right|^2-r\frac{\partial V}{\partial r} + \Delta
  f.\nonumber
  \eee
Hence the maximum of $r\frac{\partial\varphi}{\partial r}-2\varphi
  +f$ must occur on the boundary of $\partial\Omega$ which is at
  most $\max_{\partial\Omega}(-2\varphi+f)$.

\

  \begin{thm}
    For the Neumann problem with $\varphi=-\log u_1$,
    \be
    \label{equ:2.30}
    r\frac{\partial\varphi}{\partial r}-2\varphi+f\leq \max_{\partial \Omega}(f-2\varphi)
    \ee
    where $f$ is any function satisfies (\ref{equ:2.29}).
  \end{thm}

  If we normalize $u_1$ so that $u_1\leq 1$ on $\partial\Omega$ then
  $\max_{\partial\Omega}(-2\varphi)\leq 0$ and (\ref{equ:2.30})
  gives a good growth estimate of $\varphi$.

  For example, if $\frac{\partial V}{\partial r}\geq 0$, we can then
  take $f=0$ and (\ref{equ:2.30}) says that $\frac{\varphi}{r^2}$ is
  monotonic decreasing which means that $u_1$ decays like a
  Gaussian.

\

\section{Estimate of gap for more general potential}

\hspace{5mm}We shall improve the estimate that we obtained in section one.

Let $c$ be any constant greater than $\sup u$ when $u =
\frac{u_2}{u_1}$. Let $\alpha$ be a positive constant to be
determined. Then consider the function

  \be
  F = \frac{|\nabla u|^2}{(c-u)^2} + \alpha \log (c-u).
  \ee

Then
  \be
  F_i = 2(\Sigma u_j u_{ji})(c-u)^{-2} + 2 |\nabla u|^2 u_i
(c-u)^{-3} - \alpha u_i (c-u)^{-1},
  \ee

  \bee
  \Delta F &=& 2(\sum u_{ji}^2)(c-u)^{-2} + 2(\sum u_j (\Delta
  u)_j)(c-u)^{-1}\\
  &&+ 8 (\sum u_j u_{ji} u_i)(c-u)^{-3} + 2 |\nabla u|^2 \Delta
  u (c-u)^{-3} \nonumber \\
  &&+ 6 |\nabla u|^4 (c-u)^{-4} - \alpha (\Delta
  u) (c-u)^{-1} \nonumber \\
  &&- \alpha |\nabla u|^2 (c-u)^{-2} .\nonumber
  \eee

  Since $u$ satisfies the Neumann condition and $\partial \Omega$ is
  assumed to be convex, $F$ can not achieve its maximum at the boundary
  of $\Omega$ as its normal derivative would have to be positive. So
  we assume $F$ achieves its maximum in the interior of $\Omega$
  where $\nabla F = 0$.

  If $\nabla u \neq 0$ at this point, we can choose coordinate so
  that $u_1 \neq 0$ and $u_i = 0$ for $i > 1$. Then

  \be
  u_{11}(c-u)^{-1} + |\nabla u|^2 (c-u)^{-2} = \frac{\alpha}{2}.
  \ee

  Hence
  \bee
  \Delta F &\geq& 2|\nabla u|^4 (c-u)^{-4} - 2\alpha |\nabla u|^2
  (c-u)^{-2}\\
  && + \frac{\alpha ^2}{2} - 2 (\lambda _2 - \lambda _1)|\nabla u|^2
  (c-u)^{-2}\nonumber \\
  && + 4 (\inf \varphi_{ii}) |\nabla u|^2
  (c-u)^{-2}\nonumber \\
  && + 4 \alpha |\nabla u|^2
  (c-u)^{-2} - 2 |\nabla u|^4
  (c-u)^{-4}\nonumber \\
  && -2 (\lambda _2 - \lambda _1) u (c-u)^{-1}|\nabla u|^2
  (c-u)^{-2}\nonumber \\
  && +\alpha (\lambda _2 - \lambda _1) u (c-u)^{-1} - \alpha |\nabla u|^2
  (c-u)^{-2}.\nonumber
  \eee

  If we choose $\alpha$ so that

  \be
  \alpha \geq 2 (\lambda _2 - \lambda _1) - 4 \inf \varphi_{ii} + 2 (\lambda _2 - \lambda
  _1) (\sup u) (c - \sup u)^{-1},
  \ee

  \be
  \alpha > 2 (\lambda _2 - \lambda
  _1) (\sup u) (c - \sup u)^{-1}.
  \ee

  Then $\Delta F > 0$ which is not possible. Hence at $\nabla F =
  0$, $\nabla u = 0$ and we obtain

  \be
  \sup F \leq \alpha \log c.
  \ee

  If we choose $c = (1 + \varepsilon) \sup u$ with $\varepsilon >
  0$, we can choose

  \be
  \label{equ:3.9}
  \alpha =  2 (\lambda _2 - \lambda_1)(1+\varepsilon^{-1}) - 4 \inf
  \varphi_{ii}.
  \ee

  \noindent (Here we assume $\inf
  \varphi_{ii}\leq 0$, otherwise we can apply section 1.)

\

  \begin{thm}
    Choose $\alpha$ to be (\ref{equ:3.9}), then
    \be
    \frac{\left|\nabla u\right|}{c-u} \leq \sqrt{\alpha} (\log (c) -
    \log(c-u))^{\frac{1}{2}}.
    \ee
    Therefore
    \be
    \Big|\nabla \Big(\log \Big(\frac{c}{c-u}\Big)\Big)^{\frac{1}{2}}\Big| \leq
    \frac{1}{2}\sqrt{\alpha}.
    \ee
  \end{thm}

Integrating this inequality from $u = \sup u$ to $u = 0$, we find

  \be
  \sqrt{\log \Big( 1 + \frac{1}{\varepsilon} \Big)} \leq
  \frac{1}{2}\sqrt{\alpha}d(\Omega).
  \ee

  Hence
  \bes
  \alpha \geq 4 \log \Big( 1 + \frac{1}{\varepsilon} \Big) d(\Omega)^{-2}.
  \ees

  In particular
  \be
  (\lambda _2 - \lambda_1)(1 + \varepsilon^{-1}) \geq 2 \log \Big( 1 + \frac{1}{\varepsilon} \Big)
  d(\Omega)^{-2} + 2 \inf
 \varphi_{ii} .
  \ee
Hence
  \be
  \lambda _2 - \lambda_1 \geq 2 d(\Omega)^{-2} \exp[(\inf
  \varphi_{ii})d(\Omega)^2].
  \ee

\

  \begin{thm}
  \label{thm:3.2}
    Let $\Omega$ be a convex domain so that for the first
    eigenfunction $u_1$ of the operator $-\Delta + V$, the Hessian
    of $-\log u_1$ is greater than $-a$. Then the gap of the
    first eigenfunction of the operator $-\Delta + V$ is greater than
    \be
    \label{equ:3.15}
    \lambda _2 - \lambda_1 \geq 2 d(\Omega)^{-2} \exp (-a d^2 (\Omega)).
    \ee
  \end{thm}

  Note that we have estimate $a$ in section 2 already and
  (\ref{equ:3.15}) does give a gap estimate for arbitrary smooth
  potential.

  Note that Theorem \ref{thm:3.2} shows that it is possible to
  estimate $\lambda _2 - \lambda_1$ from below depending only on the
  lower bound of the Hessian of potential as long as $\Omega$ is
  convex and $d(\Omega)$ is finite. The estimate may not be optimal
  and it is possible that $d(\Omega)$ should be replaced by integral
  of some function.

\

\section{Behavior of the ground state}
 
 \hspace{5mm} It is clear from the above discussions that the behavior of the
  Hessian of the function $\varphi = -\log u_1$ is important. Since

  \be
  \Delta \varphi = |\nabla \varphi|^2 - V + \lambda_1.
  \ee
  It is clear that upper estimate of $\Delta \varphi$ can be used to
  control the growth of $\varphi$ and hence the growth of $u_1$.

  Clearly,
  \be
  \Delta(\Delta \varphi) = 2\sum \varphi^2_{ij} -2 \sum
  \varphi_j (\Delta \varphi)_j - \Delta V.
  \ee

  Let $\rho$ be a nonnegative function which varnishes on $\partial
  \Omega$, then
  \bee
  \Delta (\rho ^2 \Delta \varphi ) & = & 2 (\rho \Delta + |\nabla \rho|^2)\Delta
  \varphi + 2 \rho \nabla \rho \cdot \nabla (\Delta \varphi) \\
  && + \rho ^2 (2 \sum \varphi^2_{ij}-2 \sum
  \varphi_j (\Delta \varphi)_j - \Delta V). \nonumber
  \eee

  At the point where $\rho ^2 \Delta \varphi$ achieves its maximum,
  $\nabla( \rho^2 \Delta \varphi) = 0$ and

  \be
  \rho \nabla (\Delta \varphi) + 2 (\Delta \varphi )\nabla \rho = 0.
  \ee

  Hence
  \bee
  \Delta ( \rho^2 \Delta \varphi) & = & 2 (\rho \Delta \rho - 3|\nabla
  \rho|^2) \Delta \varphi \\
  && + 2 \rho ^2 \sum \varphi^2_{ij} - 4 \rho \Delta \varphi (\rho \cdot \nabla \varphi
  ) - \rho \Delta V.\nonumber
  \eee

  Note
  \be
  |\nabla \rho \cdot \nabla \varphi| \leq |\nabla \rho |(\sqrt{|\nabla \varphi|^2 - V + \lambda_1} + \sqrt{(V -
  \lambda_1)_+)},
  \ee
  where $(V - \lambda_1)_+$ is the positive part of $V - \lambda_1$.
  Therefore when $\rho ^2 \Delta \varphi$ achieves its maximum,
  
  \bee
  0 & \geq & 2 (\rho \Delta \rho - 3|\nabla
  \rho|^2) \rho ^2 \Delta \varphi + \frac{2}{n}(\rho ^2 \Delta
  \varphi)^2\\
  && - 4 (\rho ^2 \Delta \varphi)|\nabla \rho| (\sqrt{\rho ^2 \Delta \varphi} + \sqrt{(V -
  \lambda_1)})_+ - \rho ^4 \Delta V.\nonumber
  \eee

\

  \begin{thm}
  \label{thm:4.1}
    For any function $\rho$ vanishing at the boundary of $\Omega$,
    $\rho ^2 \Delta \varphi$ is bounded from above by $\sup (\rho \Delta \rho - 3|\nabla
  \rho|^2)$, $\sup |\nabla \rho|^2$, $\sup\rho ^2 \sqrt {(\Delta V)_+}$
  and $\sup |\nabla \rho|\sqrt{(V - \lambda_1)_+}$.
  \end{thm}

  Note that if $V$ grows at most quadratically, Theorem
  \ref{thm:4.1} shows that $\Delta \varphi$ can be bounded from
  above in terms of $(\Delta V)_+$. Since $\Delta \varphi = |\nabla \varphi|^2 - V -
  \lambda_1$, $|\varphi|$ can not grow faster than the integral of
  $\sqrt{(V - \lambda_1)_+}$ along paths tend to infinity. In
  particular for the first eigenfunction $u_1 = \exp(- \varphi_1)$,
  it cannot decay too fast.

\end{document}